\documentclass[10pt]{article}

\usepackage{graphicx}%
\usepackage{amsmath,amssymb,amsthm,upref,bm}%

\newtheorem{corollary}{Corollary}%
\newtheorem{theorem}{Theorem}%
\newtheorem{proposition}{Proposition}%
\newtheorem{definition}{Definition}%
\begin{document}

\baselineskip=4.4mm

\makeatletter

\newcommand{\MF}[1]{\mathop{#1\vrule height0.45ex width0pt}\nolimits}

\newcommand{\pd}[3][]{\mathchoice{\raise-0.5pt\hbox{$\partial$}%
\vphantom{\partial}_{\mkern-1.5mu#2}^{\mkern0.4mu#1}\mkern0.3mu}%
{\raise-0.5pt\hbox{$\partial$}%
\vphantom{\partial}_{\mkern-1.5mu#2}^{\mkern0.4mu#1}\mkern0.3mu}%
{\raise-0.5pt\hbox{$\scriptstyle\partial$}%
\vphantom{\partial}_{\mkern-1.7mu#2}^{\mkern0.1mu#1}\mkern0.1mu}%
{\raise-0.5pt\hbox{$\scriptscriptstyle\partial$}%
\vphantom{\partial}_{\mkern-1.7mu#2}^{\mkern0.1mu#1}\mkern0.1mu}#3}

\newcommand{\D}{\mathrm{d}\kern0.2pt}%
\newcommand{\ii}{\kern0.05em\mathrm{i}\kern0.05em}% 
\newcommand{\E}{\textrm{e}}%
\renewcommand{\vec}[1]{\bm{#1}}%
\newcommand{\RR}{\mathbb{R}}%
\newcommand{\Cb}{\mathbb{C}}%
\def\transp{\mathsf{T}}

\renewcommand{\Re}{\mathrm{Re}} 
\renewcommand{\Im}{\mathrm{Im}}

\def\bottomfraction{0.9}

\title{\bf On a cylinder freely floating \\ in oblique waves}

\author{Nikolay Kuznetsov}

\date{}

\maketitle

\begin{center}
Laboratory for Mathematical Modelling of Wave Phenomena, \\ Institute for Problems
in Mechanical Engineering, Russian Academy of Sciences, \\ V.O., Bol'shoy pr. 61,
St. Petersburg 199178, Russian Federation \\ E-mail address:
nikolay.g.kuznetsov@gmail.com
\end{center}

\begin{abstract}
The coupled motion is investigated for a mechanical system consisting of water and a
body freely floating in it. Water occupies either a half-space or a layer of
constant depth into which an infinitely long surface-piercing cylinder is immersed,
thus allowing us to study the so-called oblique waves. Under the assumption that the
motion is of small ampli\-tude near equilibrium and describes time-harmonic
oscillations, the phe\-nomenon's linear setting reduces to a spectral problem with
the radian frequency as the spectral parameter. If the radiation condition holds,
then the total energy is finite and the equipartition of kinetic and poten\-tial
energy holds for the whole system. On this basis, it is proved that no wave modes
are trapped under some restrictions on their frequencies; in the case when a
symmetric cylinder has two immersed parts restrictions are imposed on the type of
mode as well.
\end{abstract}

\setcounter{equation}{0}

\section{Introduction}

This paper continues the author's studies dealing with the motion of a mechanical
system consisting of a water layer of constant depth and a rigid body freely
floating in it; the initial note \cite{NGK10} was written on the occasion of
V.\kern2pt M. Babich's 80th birthday, who recently turned 90. It was F.~John
\cite{John1}, who proposed the linear problem describing the coupled motion of water
bounded from above by the atmosphere and a partially immersed body. The latter
floats freely according to Archimedes' law being unaffected by all external forces
(for example due to constraints on its motion) except for gravity. The motion of
water is assumed to be irrota\-tional (viscosity is neglected as well as the surface
tension), whereas the motion of the whole system is supposed to be of small
amplitude near equilibrium; this allows us to use a linear model.

In the framework of the linear theory of water waves, two- and three-dimensional
statements of the problem are possible; the one considered here is two-dimensional
which is essential for formulation of restrictions to be imposed when investigating
the question of uniqueness. The latter is of paramount importance (it is the first
one on the list of important open problems in the survey \cite{U}), because there
are examples of non-uniqueness; see, for example,~\cite{NGK10}. For this reason
restrictions on the frequency range and on the body's shape are required. It should
be noted that such restrictions are unnecessary for similar acoustical problems;
see, for example, the monograph \cite{CK}.

The original John's formulation of the floating body problem is rather cumbersome
because he did not use matrices to express the equations of body's motion (the
matrix form of these equations described below demonstrates their simple structure;
see also \cite{KM}, \cite{KM1} and \cite{KM2}). Anyway, the problem was neglected by
researches during 60 years after publication of the article \cite{John}, in which
the uniqueness theorem was proved for the three-dimensional problem under the
assumptions that the so-called John condition holds for the body (it is described
below for two-dimensional geometry) and the frequency of oscillations is
sufficiently large.

On the other hand, the problem of time-harmonic oscillations of water in the
presence of fixed rigid bodies attracted much attention in the second half of the
20th century; initially, the case of a single body had been investigated in the
article \cite{John} mentioned above. Numerous results about this problem were
presented in detail in the summarising monograph \cite{KMV}, which also contains
extensive literature. In particular, the first non-uniqueness example due to
M.~McIver \cite{MMI} was generalized in \cite{KMV}. Namely, examples of non-trivial
solutions were constructed for the two-dimensional homogeneous problem with an
arbitrary number of fixed surface-piercing bodies (only two bodies were considered
in \cite{MMI}); the so-called inverse procedure was applied for this purpose. The
presence of multiple bodies violates the original John's condition according to
which only a single surface-piercing body is admissible in the two-dimensional
problem. From the hydrodynamic viewpoint, these non-trivial solutions describe
trapped modes, that is, free oscillations of water having final energy. These
results were developed further in the article \cite{4}, in which the uniqueness
theorem due to John and McIver's non-uniqueness example were extended to the case of
infinitely long, sur\-face-piercing cylinders in oblique waves. The problem for
fixed cylinders considered in \cite{4} is a particular case of that formulated below
in \S~2. Another approach to the question of uniqueness in the problem with fixed
cylinders was proposed in \cite{K}.

\vspace{-1mm}

\section{Statement of the problem}

Let an infinitely long surface-piercing cylinder of uniform cross-section float
freely in water which is either infinitely deep or bounded from below by a
horizontal rigid bottom. The Cartesian coordinate system $(x,y)$ is chosen in a
plane orthogonal to the cylinder's generators (directed along the $z$-axis), so
that the $y$-axis is directed upwards, whereas the mean free surface of the water
intersects this plane along the $x$-axis. Thus, the cross-section $W$ of the water
domain is a subset of $\RR^2_- = \{ x \in \RR, \, y<0 \}$. Let $\widehat{B}$ denote
the bounded two-dimensional domain whose closure is the cross-section of the
cylinder's equilibrium position; we suppose that $\widehat{B} \setminus
\overline{\RR^2_-}$\,---\,the part of the body located above the water
surface\,---\,is a nonempty domain, whereas the immersed part $B = \widehat{B} \cap
\RR^2_-$ is the union of a finite number of domains. Thus, $D = \widehat{B} \cap
\partial \RR^2_-$ consists of the same number of nonempty intervals of the $x$-axis;
see Fig.~1 for the case of two immersed parts, where
\[ D = \{ x \in (-a, -b) \cup (b, a) , \, y=0 \} .
\]
Notice that $W = \RR^2_- \setminus \skew2\overline{B}$, when the water has infinite
depth,~or
\[ W = \{ x \in \RR, \, -h<y<0 \} \setminus \skew2\overline{B} , \quad \mbox{where} \ 
h > b_0 = \sup_{(x,y) \in B} |y| ,
\]
when the water has the constant finite depth $h$ (see Fig.~1). The cross-section of
the bottom is denoted by $H = \{ x \in \RR , \, y = - h \}$ in the last case.
Furthermore, $W$ is assumed to be a Lipschitz domain, and so the unit normal
$\mathbf{n}$ pointing to the exterior of $W$ is defined almost everywhere on
$\partial W$. Finally, we denote by $S = \partial \widehat{B} \cap \RR^2_-$ the
wetted contour (the number of its components is equal to the number of immersed
domains), whereas $F = \partial \RR^2_- \setminus \skew2 \overline {D}$ is the free
surface at rest.

\begin{figure}
\begin{center}
\vspace{4mm}
\unitlength=0.6mm

\special{em:linewidth 0.4pt}

\linethickness{0.4pt}

\begin{picture}(163.00,62.00)(-3,90)

\put(28.0,130.00){\line(0,1){15.00}}

\put(65.0,130.00){\line(0,1){8.00}}

\put(89.0,130.00){\line(0,1){8.00}}

\put(126.0,130.00){\line(0,1){15.00}}

\put(28.0,145.00){\line(1,0){98.00}}

\put(65.0,138.00){\line(1,0){24.00}}

\put(100.0,120.00){\line(1,0){28.00}}

\put(26.0,120.00){\line(1,0){28.00}}

\put(-3,130.00){\vector(1,0){160.00}}

\put(77.00,130.00){\vector(0,1){24.00}}

\multiput(65,130)(0,-7.2){5}{\line(0,-1){6}}

\multiput(89,130)(0,-7.2){5}{\line(0,-1){6}}

\bezier{464}(28.00,130.00)(15.00,83.00)(65.00,130.00)

\multiput(28,130)(-6,-6){5}{\line(-1,-1){5.1}}

\multiput(126,130)(6,-6){5}{\line(1,-1){5.1}}

\bezier{464}(89.00,130.00)(139.00,83.00)(126.00,130.00)

\put(-3,95){\line(1,0){160}}

\put(10.00,134.00){\makebox(0,0)[cc]{{\small $F_{\infty}$}}}

\put(144.00,134.00){\makebox(0,0)[cc]{{\small $F_{\infty}$}}}

\put(81.00,152.00){\makebox(0,0)[cc]{{\small $y$}}}

\put(92.8,134.00){\makebox(0,0)[cc]{{\small $+b$}}}

\put(130.00,134.00){\makebox(0,0)[cc]{{\small $+a$}}}

\put(24.00,134.00){\makebox(0,0)[cc]{{\small $-a$}}}

\put(60.8,134.00){\makebox(0,0)[cc]{{\small $-b$}}}

\put(5.00,120.00){\makebox(0,0)[cc]{{\small $W_{\infty}$}}}

\put(37.00,124.00){\makebox(0,0)[cc]{{\small $B_-$}}}

\put(36.00,116.00){\makebox(0,0)[cc]{{\small ballast}}}

\put(118.00,116.00){\makebox(0,0)[cc]{{\small ballast}}}

\put(106.00,131.40){\makebox(0,0)[cc]{{\large $\widehat{B}$}}}

\put(117.00,124.00){\makebox(0,0)[cc]{{\small $B_+$}}}

\put(46.00,108.00){\makebox(0,0)[cc]{{\small $S_-$}}}

\put(110.00,108.00){\makebox(0,0)[cc]{{\small $S_+$}}}

\put(152.00,120.00){\makebox(0,0)[cc]{{\small $W_{\infty}$}}}

\put(155.00,134.00){\makebox(0,0)[cc]{{\small $x$}}}

\put(152.00,98.00){\makebox(0,0)[cc]{{\small $-h$}}}

\put(77.00,111.00){\makebox(0,0)[cc]{{\small $W_0$}}}

\put(82.00,134.00){\makebox(0,0)[cc]{{\small $F_0$}}}

\bezier{28}(128,122)(132.70,124)(134.70,130.00)

\put(137.50,123.50){\makebox(0,0)[cc]{{\small $\beta$}}}

\end{picture}

\end{center}
\vspace{-6mm} \caption{A definition sketch of the cylinder's cross-section
$\widehat{B}$ with two immersed parts; they are denoted by $B_-$ and $B_+$, and
their wetted boundaries are $S_-$ and $S_+$ respectively. The cross-section $F$ of
the free surface of the water consists of three parts; two of them lying on the
$x$-axis outside $|x| > a$ are denoted by $F_\infty$ and the third one $F_0$ is
between $x = -b$ and $x = +b$. Furthermore, $W = W_0 \cup W_\infty$, where $W_0$ is
the part of the water domain located in the vertical strip under $F_0$ and
$W_\infty$ is its complement. The equation of the horizontal bottom is $y = -h$.}
\vspace{-4mm}
\end{figure}
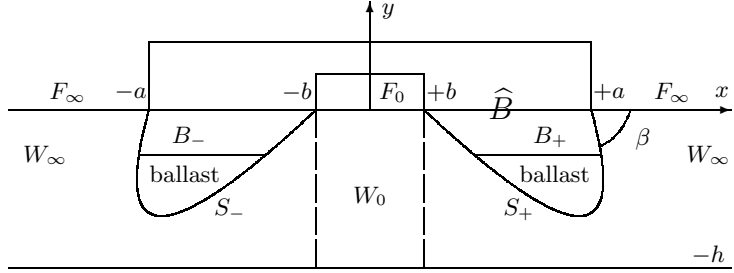

To describe the small-amplitude coupled motion of the water/cylinder system it is
standard to apply the linear setting, in which case the first-order approximation of
unknowns is used. These are the velocity potential $\MF{\Phi} (x,y,z;t)$ and the
vector column $\MF{\mathbf{q}} (t)$ describing the motion of the cylinder; its three
components are as follows:

\vspace{1mm}

\noindent $\bullet$ $q_1$ and $q_2$ are the displacements of the centre of mass in
the horizontal and vertical directions respectively from its rest position $\bigl(
x^{(0)}, y^{(0)} \bigr)$;

\noindent $\bullet$ $q_3$ is the angle of rotation about the axis that goes through
the centre of mass orthogonally to the $(x,y)$-plane (the angle is measured from the
$x$- to the $y$-axis).

\vspace{1mm}

\noindent We omit relations governing the time-dependent behaviour (details can be
found in \cite{NGK10}) and turn to the description of time-harmonic oscillations of
the coupled water/cylinder system in the presence of oblique waves for which purpose
the following ansatz
\begin{equation}
 \bigl( \MF{\Phi} (x,y,z; t), \mathbf{q} (t) \bigr) = \Re \bigl\{ \bigl( \E^{\ii(k z 
 - \omega t)} \MF{\varphi} (x,y), \ii \E^{-\omega t} \mathbf{z} \bigr) \bigr\} 
 \label{eq:ansatz}
\end{equation}
is applied. Here $\omega > 0$ is the radian frequency of oscillations, to which the
wavenumber $\nu = \omega^2 / g$ corresponds; $g > 0$ is the acceleration due to
gravity that acts in the direction opposite to the $y$-axis. Furthermore, $k \in [0,
\nu)$ is the prescribed wavenumber component that the wave train has parallel to the
generators of the cylinder; $\varphi \in \MF{H^1_{loc}}(W)$ is a complex-valued
function and $\mathbf{z} \in \Cb^3$. Thus, $k / \nu$ is the sine of the angle
between the wave crests and the plane normal to the generators; waves are called
oblique when $k > 0$.

To be specific, we consider the case of infinitely deep $W$ first. In the absence
of incident waves, we obtain the following problem for $(\varphi, \mathbf{z})$:
\begin{gather}
 (\nabla^2 - k^2) \varphi = 0 \quad \mbox{in} \ W ,
 \label{eq:1}\\
 \pd{y} \varphi - \nu \varphi = 0 \quad \mbox{on} \ F, 
 \label{eq:2}\\
 \pd{\mathbf{n}} \varphi = \omega \, \mathbf{N}^\transp \mathbf{z} \ \Big(\!\! = 
 \omega \sum_1^3 N_j z_j \Big) \quad \mbox{on} \ S , \label{eq:4}\\
 \nabla \varphi \to 0 \quad \mbox{as} \ y \to -\infty ,
 \label{eq:3}\\
 \omega^2 \mathbf{E} \mathbf{z} = - \omega \int_{S} \varphi \, \mathbf{N} \,
 \D{}s + g \, \mathbf{K} \mathbf{z} . \label{eq:5}
\end{gather}
Here $\nabla = (\pd{x}, \pd{y})$ is the spatial gradient, whereas $\mathbf{N} =
(N_1, N_2, N_3)^\transp$ (the operation $^\transp$ transforms a vector row into a
vector column and vice versa), where $(N_1, N_2)^\transp = \mathbf{n}$, $N_3 =
\left( x - x^{(0)}, y - y^{(0)} \right)^\transp \times \mathbf{n}$ and $\times$
stands for the vector product. In the equations of the body motion \eqref{eq:5}, the
$3\!\times\!3$ matrices are as follows:
\begin{equation}
 \MF{\mathbf{E}} =
\begin{pmatrix}
 I^M & 0 & 0 \\
 0 & I^M & 0 \\
 0 & 0 & I^M_2
\end{pmatrix} \quad {\rm and} \quad \MF{\mathbf{K}} =
\begin{pmatrix}
 0 & 0 & 0 \\
 0 & I^D & I^D_x \\
 0 & I^D_x & I^D_{xx} + I^B_y
\end{pmatrix} .
\label{eq:EK}
\end{equation}
The positive elements of the mass/inertia matrix $\MF{\mathbf{E}}$ are as follows:
\begin{gather*}
I^M = \rho_0^{-1} \int_{\widehat{B}} \MF{\rho}(x,y) \, \D x \D y \, , \\ I^M_2 =
\rho_0^{-1} \int_{\widehat{B}} \MF{\rho}(x,y) \Big[ \left( x - x^{(0)} \right)^2 +
\left( y - y^{(0)} \right)^2 \Big] \D x \D y \, .
\end{gather*}
Here $\MF{\rho} (x,y) \geq 0$ is the density distribution within the body and
$\rho_0 > 0$ is the constant density of water. On the right-hand side of relation
\eqref{eq:5}, we have forces and their moments: the first term is due to the
hydrodynamic pressure and the second one is related to the buoyancy (see, for
example, \cite{John1}). The non-zero elements of the matrix $\mathbf{K}$ are
\begin{gather*}
I^D = \int_D \D x > 0, \quad I^D_x = \int_D \big( x - x^{(0)} \big) \D x , \\
I^D_{xx} = \int_D \big( x - x^{(0)} \big)^2 \D x > 0 , \quad I^B_y = \int_B \big( y
- y^{(0)} \big) \D x \, \D y .
\end{gather*}
It should be noted that the matrix $\mathbf{K}$ is symmetric.

In relations \eqref{eq:2}, \eqref{eq:4} and \eqref{eq:5}, $\omega$ is a spectral
parameter which is sought together with the eigenvector $(\varphi,\mathbf{z})$.
Since $W$ is a Lipschitz domain and $\varphi \in \MF{H^1_{loc}}(W)$, relations
\eqref{eq:1}--\eqref{eq:4} are understood in the sense of the following integral
identity:
\begin{equation}
 \int_{W} \nabla \varphi \nabla \psi \,\D x \D{}y = \nu 
 \int_{F} \varphi \, \psi \, \D x + \omega \int_{S} \psi \, 
 \mathbf{N}^\transp \mathbf{z} \, \D{}s .
\label{eq:intid}
\end{equation}
It must hold for an arbitrary smooth $\psi$ having a compact support in~$\overline
W$.

Together with \eqref{eq:3} the following condition
\begin{equation}
\int_{W\cap\{|x|=b\}} \bigl| \pd{|x|} \varphi - \ii \ell \varphi \bigr|^2\, \D{}s =
\MF{o}(1) \quad \mbox{as} \ b \to \infty , \ \ \ell = (\nu^2 - k^2)^{1/2} ,
\label{eq:6}
\end{equation}
specifies the behaviour of $\varphi$ at infinity; \eqref{eq:3} means that the
velocity field decays with depth, whereas \eqref{eq:6} yields that the potential
given by formula \eqref{eq:ansatz} describes outgoing waves. This radiation
condition is similar to that used in \cite{John}, where the problem was considered
for two-dimensional water waves in the presence of a fixed obstacle.

The relations listed above must be augmented by the following subsi\-diary
conditions concerning the equilibrium position:

\vspace{1mm}

\noindent $\bullet$ Archimedes' law, $I^M = \int_B \D x \D{}y$ (the mass of the
displaced liquid is equal to that of the body);

\noindent $\bullet$ $\int_B \bigl(x - x^{(0)} \bigr) \, \D x \D{}y = 0$ (the centre
of buoyancy lies on the same vertical line as the centre of mass);

\noindent $\bullet$ The matrix $\mathbf{K}$ is positive semi-definite; moreover, the
$2 \times 2$ matrix $\mathbf{K}'$ that stands in the lower right corner of
$\mathbf{K}$ is positive definite (see \cite{John1}).

\vspace{1mm}

\noindent The last of these requirements yields stability of the body equilibrium
position, which follows from the results formulated, for example, in the paper
\cite[\S\,2.4]{John1}. The stability is understood in the classical sense, that is,
an instantaneous infinitesimal disturbance causes the position changes which remain
infinitesimal, except for purely horizontal drift, for all subsequent times.

In conclusion of this section, we note that relations \eqref{eq:3} and \eqref{eq:6}
must be amended in the case when $W$ has finite depth. Namely, the no flow condition
\begin{equation}
\pd{y} \varphi = 0 \quad \mbox{on} \ H \label{eq:H}
\end{equation}
replaces \eqref{eq:3}, whereas $\ell$ must be changed to $\ell_0$ in \eqref{eq:6},
where $\ell_0$ is the unique positive root of $\ell_0 \MF{\tanh}(\ell_0 h) = \ell$.

\section{Equipartition of energy}

It is known (see, for example, \cite[\S\,2.2.1]{KMV}), that a potential, satisfying
relations \eqref{eq:1}, \eqref{eq:2}, \eqref{eq:3} and \eqref{eq:6}, has an
asymptotic representation at infinity of the same type as Green's function. Namely,
if $W$ has infinite depth, then
\begin{eqnarray}
&& \ \ \ \ \ \ \ \ \MF{\varphi} (x, y) = \MF{A}_\pm (y) \, \E^{\ii \ell |x|} +
\MF{r}_\pm (x, y) , \nonumber \\ && |r_\pm|^2, \, |\nabla r_\pm| = \MF{O} \bigl(
 [x^2 + y^2]^{-1} \bigr) \ \mbox{as} \ x^2 + y^2 \to \infty , \label{eq:uas}
\end{eqnarray}
and the following equality holds
\begin{equation}
 \ell \int_{-\infty}^0 \left( |\MF{A}_+ (y)|^2 + |\MF{A}_- (y)|^2 \right) \D{} y =
 -\Im \int_S \overline{\varphi}\,\pd{\mathbf{n}} \varphi \, \D{}s .
\label{eq:ener}
\end{equation}

Assuming that $( \MF{\varphi} , \mathbf{z} )$ is a solution of the problem
\eqref{eq:1}--\eqref{eq:5} and \eqref{eq:6}, we rearrange the last formula using
the coupling conditions \eqref{eq:4} and \eqref{eq:5}. First, transposing the
complex conjugate of equation \eqref{eq:5}, we obtain
\[ \omega^2 \left( \mathbf{E} \overline{\mathbf{z}} \right)^\transp = - \omega
 \int_{S} \overline{\varphi} \, \mathbf{N}^\transp \, \D{}s + g \left( \mathbf{K}
 \overline{\mathbf{z}} \right)^\transp .
\]
This relation and condition \eqref{eq:4} yield that the inner product of both sides
with $\mathbf{z}$ can be written in the form
\begin{equation}
\omega^2 \, \overline{\mathbf{z}}^\transp \mathbf{E} \mathbf{z} - g
\overline{\mathbf{z}}^\transp \mathbf{K} \mathbf{z} = - \int_{S} \overline{\varphi}
\, \pd{\mathbf{n}} \varphi \,\D{}s .
\label{eq:transp}
\end{equation}
Second, substituting this equality into \eqref{eq:ener}, we obtain
\begin{equation}
 \ell \int_{-\infty}^0 \left( |\MF{A}_+ (y)|^2 + |\MF{A}_- (y)|^2 \right) \D{} y =
 \Im \Bigl\{ \omega^2 \, \overline{\mathbf{z}}^\transp \mathbf{E} \mathbf{z} - g 
 \overline{\mathbf{z}}^\transp \mathbf{K} \mathbf{z} \Bigr\} . \label{eq:ener'}
\end{equation}
In the same way as in \cite{KM1}, this yields the following assertion about the
kinetic and potential energy of the water motion.

\begin{proposition}\label{propos:1}
Let\/ $( \MF{\varphi} , \mathbf{z} )$ be a solution of problem\/
\eqref{eq:1}--\eqref{eq:5} and\/ \eqref{eq:6}, then
\begin{equation}
 \int_W \big( |\nabla \varphi|^2 + k^2 |\varphi|^2 \big) \,\D{}x \, \D{}y < \infty
 \quad \mbox{and} \quad
 \nu \int_F |\varphi|^2 \, \D{}x < \infty \, , \label{eq:finenerg}
\end{equation}
that is, $\varphi \in H^1 (W)$. Moreover, the following equality holds:
\begin{equation}
 \int_W \big( |\nabla \varphi|^2 + k^2 |\varphi|^2 \big) \,\D{}x \, \D{}y + \omega^2
 \overline{\mathbf{z}}^\transp \mathbf{E} \mathbf{z} = \nu \int_F |\varphi|^2 \, \D{}x
 + g \, \overline{\mathbf{z}}^\transp \mathbf{K} \mathbf{z} . \label{eq:lagrange}
\end{equation}
\end{proposition}

Here the kinetic energy of the water/body system stands on the left-hand side,
whereas we have the potential energy of this coupled motion on the right-hand side.
Thus the last formula generalises the energy equiparti\-tion equality valid when a
fixed body is immersed into water. Indeed, $\mathbf{z}$ vanishes for such a body,
and \eqref{eq:lagrange} turns into the well-known equality (see, for example,
formula (4.99) in \cite{KMV}).

Proposition 1 shows that if $( \varphi, \mathbf{z} )$ is a solution of problem
\eqref{eq:1}--\eqref{eq:5} and \eqref{eq:6} with complex-valued components, then its
real and imaginary parts separately satisfy this problem. This allows us to consider
$( \varphi, \mathbf{z} )$ as an element of the real product space $H^1 (W) \times
\RR^3$ in what follows.

\begin{definition}\rm
Let the subsidiary conditions concerning the equilibrium position (see \S~2) hold
for the freely floating body $\widehat{B}$. A non-trivial real solution $( \varphi,
\mathbf{z} ) \in H^1 (W) \times \RR^3$ of problem \eqref{eq:intid} and \eqref{eq:5}
is called a {\it mode trapped}\/ by this body, whereas the corresponding value of
$\omega$ is referred to as a {\it trapping frequency}.
\end{definition}

In the case of finite depth, the remainder in formula \eqref{eq:uas} has the
following behaviour uniformly in $y \in [-h, 0]$:
\[ |r_\pm (x, y)| , \ \ |\nabla r_\pm (x, y)| = \MF{O} \bigl( |x|^{-1} \bigr) \quad 
\mbox{as} \ |x| \to \infty \, ,
\]
whereas formula \eqref{eq:ener} holds with $\ell$ changed to $\ell_0$. Therefore,
Proposition~1 is true for problem \eqref{eq:1}--\eqref{eq:5} and \eqref{eq:6} with
condition \eqref{eq:3} replaced by \eqref{eq:H} and $\ell$ changed to $\ell_0$ in
\eqref{eq:6}. Definition~1 remains unchanged for the finite depth case.

\section{On the absence of trapped modes}

In order to determine conditions on the domain $\widehat{B}$ and on the frequency
$\omega$ guaranteeing the absence of non-trivial solutions $( \varphi, \mathbf{z})
\in H^1 (W) \times \RR^3$, we write \eqref{eq:lagrange} as follows:
\begin{equation}
\mathbf{z}^\transp ( \omega^2 \mathbf{E} - g \, \mathbf{K} ) \mathbf{z}  = \nu
\int_F |\varphi|^2 \, \D{}x - \int_W |\nabla \varphi|^2\,\D{}x \D{}y .
\label{eq:eq}
\end{equation}
It is clear that the left-hand side is non-negative for a nonzero $\mathbf{z}$
provided $\omega^2 \geq \lambda_0$; the latter is the largest $\lambda$ satisfying
$\det( \lambda \mathbf{E} - g \mathbf{K} ) = 0$. Let the equilibrium position be
stable for the cylinder, whose cross-section is $\widehat{B}$ (see \S~2 for
conditions of stability), it is convenient to say that property $\Omega$ holds for
$\omega$ if $\omega^2 \geq \lambda_0$. Thus, we arrive at the following.

\begin{proposition}\label{propos.2}
Let the cylinder's cross-section be $\widehat{B}$ and let property $\Omega$ hold for
$\omega$. If the inequality
\begin{equation}
\nu \int_F |\varphi|^2 \, \D{}x < \int_W \big( |\nabla \varphi|^2 + k^2 |\varphi|^2
\big) \,\D{}x \D{}y \, ,
\label{eq:ineq}
\end{equation}
holds for every non-trivial $\varphi \in H^1 (W)$ satisfying \eqref{eq:intid} and
\eqref{eq:5}, then $\omega$ is not a trapping frequency.
\end{proposition}

Indeed, for a non-trivial $( \varphi, \mathbf{z} )$ inequality \eqref{eq:ineq} is
incompatible with property $\Omega$.

\subsection{Cylinders with single immersed part}

For the sake of simplicity we suppose that $W$ has infinite depth. We recall that
John's condition requires a simply connected domain $\widehat{B} \cap \RR^2_-$ to
belong to the strip lying between the verticals going through the points, where the
contour $\partial \widehat{B}$ intersects the $x$-axis.

\begin{theorem}\label{theor.1}
Let $W$ have infinite depth and let the domain $\widehat{B} \cap \RR^2_-$ be simply
connected and satisfy John's condition. Let also the subsidiary conditions
concerning the equilibrium position (see \S~2) hold. Then inequality
\eqref{eq:ineq} holds for a non-trivial solution $(\varphi, \mathbf{z}) \in H^1 (W)
\times \RR^3$ of problem \eqref{eq:intid} and \eqref{eq:5}.
\end{theorem}

\noindent {\it Proof.} Since $\widehat{B} \cap \RR^2_-$ is simply connected, $D$
consists of a single interval, and so the intersection of the free surface and the
$(x, y)$-plane is $F = F_\infty$ (see Fig.~1); moreover, $F = F_+ \cup F_-$, where
$F_+$ $(F_-)$ is the ray lying on the $x$-axis to the right (left) of $D$. Let us
prove the inequality
\begin{equation}
\nu \int_{F_\pm} |\varphi|^2 \, \D{}x < \int_{W_\pm} \big( |\nabla \varphi|^2 + k^2
|\varphi|^2 \big) \,\D{}x \D{}y \, ,
\label{eq:ineq+}
\end{equation}
where $W_\pm \! \subset \! W$ is the subdomain lying strictly under $F_\pm$; in view
of John's condition these subdomains are well defined.

Following John let us define
\begin{equation}
a^{(\pm)} (x) = \int_{-\infty}^0 \MF{\varphi (x,y)} \E^{\nu y}\, \D y \quad
\mbox{on} \ F_\pm \, . \label{14.6}
\end{equation}
Differentiating this function twice and using equation \eqref{eq:1}, we obtain
\[  a_{xx}^{(\pm)} = k^2 a^{(\pm)} - \int_{-\infty}^0 \varphi_{yy} (x,y) \,
\E^{\nu y}\, \D y \, .
\]
Integrating by parts twice and taking into account conditions \eqref{eq:2} and
\eqref{eq:3}, we conclude that $a_{xx}^{(\pm)} + \ell^2 a^{(\pm)} = 0$ on $F_\pm$
because $k^2 - \nu^2 = - \ell^2 $. 

Since $\varphi \in H^1 (W)$, we have that $\lim_{|x| \to \infty} \MF{a^{(\pm)}}(x) =
0$, and so $a^{(\pm)} \equiv 0$ on $F_{\pm}$. Now, integrating by parts in
\eqref{14.6}, we see that
\[ \varphi (x,0) = \int_{-\infty}^0 \MF{\varphi_y} (x,y) \, \E^{\nu y} \, \D y \quad
\textrm{for} \ (x,0) \in F_\pm \, .
\]
Squaring both sides, applying the Schwarz inequality to the integral and integrating
over $F_\pm$, we obtain
\[ \nu \int_{F_\pm} |\phi^{(\pm)}(x,0)|^2\, \D x  \leq \frac{1}{2} \int_{W_\pm}
|\MF{\varphi_y}|^2\, \D x\/ \D y < \int_{W_\pm} |\nabla \varphi|^2 \, \D x\/ \D y
\, ;
\]
here the coefficient $1/2$ results from integration of $\E^{2 \nu y}$. The last
inequality is stronger than \eqref{eq:ineq+}, and so it implies \eqref{eq:ineq}
because $F = F_+ \cup F_-$. $\square$

\vspace{1mm}

Propositions 1 and 2, and Theorem 1 yield.

\begin{corollary}\label{corol.1}
Let $W$ have infinite depth and let the domain $\widehat{B} \cap \RR^2_-$ be simply
connected and satisfy John's condition. Let also the subsidiary conditions
concerning the equilibrium position (see \S~2) hold. If property $\Omega$ holds for
$\omega$, then problem \eqref{eq:1}--\eqref{eq:5} and\/ \eqref{eq:6} has only a
trivial solution for this $\omega$.
\end{corollary}

Analogous assertion is true when the water domain has finite depth.

\subsection{Cylinders with two immersed parts}

In the case of a cylinder with two immersed parts (see Fig.~1), we are going to
apply the trick based on John's condition which was used in \S~4.1. This is possible
if the cylinder's cross-section is symmetric about a vertical axis (the $y$-axis in
Fig.~1), the density distribution $\rho$ is also symmetric within the cylinder and
each of the domains $B_+$ and $B_-$ is simply connected and satisfies John's
condition, that is, $\beta \geq \pi / 2$ (see Fig.~1). Let also the subsidiary
conditions concerning the equilibrium position (see \S~2) hold. Moreover, we
restrict the cylinder's admissible motions to sway, heave and roll, and consider
separately the classes of velocity potentials consisting of odd and even functions
of $x$.

\subsubsection{Sway motion}

In this case $\mathbf{z} = (z_1, 0, 0)^\transp$, $z_1 \in \RR$, and so the boundary
condition \eqref{eq:4} takes the form
\begin{equation}
\pd{\mathbf{n}} \varphi = \omega N_1 z_1 \quad \mbox{on} \ S \, , \label{eq:4'}
\end{equation}
and system \eqref{eq:5} splits. Its first equation is as follows
\begin{equation}
\omega I^M z_1 = - \int_{S} \varphi \, N_1 \, \D{}s \, , \label{eq:5'}
\end{equation}
whereas the second and third equations turn into the orthogonality condi\-tions:
\begin{equation}
\int_{S} \varphi \, N_j \, \D{}s = 0 \, , \quad j = 2, 3 . \label{eq:5''}
\end{equation}
Thus, we have to determine frequency intervals for which the equality
\begin{equation}
I^M (\omega z_1)^2 = \nu \int_F |\varphi|^2 \, \D{}x - \int_W \big( |\nabla
\varphi|^2 + k^2 |\varphi|^2 \big) \,\D{}x \D{}y \label{eq:eq+}
\end{equation}
obtained from \eqref{eq:eq} cannot hold for a non-trivial $\varphi$.

Let us assume that the domains $B_+$ and $B_-$ are simply connected and each
satisfies John's condition, that is, $\beta \geq \pi / 2$ (see Fig.~1). Repeating
literally considerations used in \S~4.1, we arrive at the inequality:
\begin{equation}
\nu \int_{F_\infty} |\varphi|^2 \, \D{}x < \int_{W_\infty} \big( |\nabla \varphi|^2
+ k^2 |\varphi|^2 \big) \,\D{}x \D{}y \, . \label{eq:ineq'}
\end{equation}
It remains to estimate $\nu \int_{F_0} |\varphi|^2 \, \D{}x$, and for this purpose
we consider the function
\begin{equation}
a (x) = \int_{-\infty}^0 \MF{\varphi (x,y)} \E^{\nu y}\, \D y \, , \label{14.6-'}
\end{equation}
which, in view of John's condition, is well defined for $x \in F_0$. In the same
way as in \S~4.1, we obtain that it satisfies the equation $a_{xx} + \ell^2 a = 0$
on~$F_0$.

Assuming that $\varphi (x,y)$ is odd in $x$, we get that $a (x) = C \sin \ell x$,
which we substitute into \eqref{14.6-'} and differentiate. After squaring both
obtained equalities, we apply the Schwarz inequality to the integrals. In this way,
we find that for every $(x, 0) \in F_0$:
\begin{gather}
2 \, \nu \, C^2 \sin^2 \ell x \leq \int_{-\infty}^0 | \MF{\varphi} (x,y) |^2 \, \D y
\, , \label{8.11-'} \\ 2 \, \nu \, \ell^2 \, C^2 \cos^2 \ell x \leq \int_{-\infty}^0
| \MF{\varphi_x} (x,y) |^2 \, \D y \, . \label{8.11-''}
\end{gather}

Furthermore, integration by parts in \eqref{14.6-'} yields
\[ \varphi (x,0) = - \nu C \sin \ell x + \int_{-\infty}^0 \MF{\varphi_y}
(x,y) \, \E^{\nu y} \, \D y 
\]
for every $(x, 0) \in F_0$, and so
\begin{equation}
\nu |\varphi (x,0)|^2 \leq 2 \nu^3 C^2 \sin^2 \ell x + \int_{-\infty}^0
|\MF{\varphi_y} (x,y)|^2 \, \D y \, .
\label{14.10}
\end{equation}
Integrating this over $F_{0}$ and using \eqref{8.11-'}, we obtain
\begin{multline}
\nu \int_{F_0} |\varphi (x,0)|^2\, \D x \leq 2 \nu (k^2 + \ell^2) C^2 \int_{F_0}
\sin^2 \ell x \, \D x + \int_{W_0} |\varphi_y|^2\, \D x \D y \\ \leq k^2 \int_{W_0}
|\varphi|^2 \, \D x \D y + \int_{W_0} |\varphi_y|^2 \, \D x \D y + 2 \nu \ell^2 C^2
\int_{F_0} \sin^2 \ell x \, \D x \, . \label{new}
\end{multline}

Let us say that property $\Omega_-$ holds for $\omega$ if the inequalities
\[ \pi m \leq \ell \, b \leq \pi (2 m + 1) / 2 
\]
are valid for some $m = 0,1,\dots$; we recall that $\ell = (\nu^2 - k^2)^{1/2}$,
$\nu = \omega^2 / g$ and $2 b$ is the spacing between $B_+$ and $B_-$.

Since property $\Omega_-$ is equivalent to the inequality
\begin{equation*}
\int_0^b \sin^2 \ell x \, \D x \leq \int_0^b \cos^2 \ell x \, \D x \, ,
\label{14.12}
\end{equation*}
we can estimate the last term in \eqref{new} with the help of \eqref{8.11-''}
provided this property holds. In this way, we arrive at the estimate
\[ \nu \int_{F_0} |\varphi|^2 \, \D{}x \leq \int_{W_0} \big( |\nabla \varphi|^2
+ k^2 |\varphi| ^2 \big) \,\D{}x \D{}y \, ,
\]
which combined with \eqref{eq:ineq'} leads to a contradiction with equality
\eqref{eq:eq+}, unless $(\varphi, \mathbf{z})$ is trivial. This completes the proof
of the following.

\begin{proposition}\label{propos.3}
Let $W$ have infinite depth and let the domain $\widehat{B}$ be sym\-metric about
the $y$-axis and such that $\widehat{B} \cap \RR^2_-$ is the union of two simply
connected domains each satisfying John's condition. Let also the subsidiary
conditions concerning the equilibrium position (see \S~2) hold.

If property $\Omega_-$ holds for $\omega$, then problem \eqref{eq:1}--\eqref{eq:5}
and\/ \eqref{eq:6} has only a trivial solution $(\varphi, \mathbf{z})$ for this
$\omega$ provided $\varphi$ is odd in $x$ and $\mathbf{z} = (z_1, 0, 0)^\transp$.
\end{proposition}

Let the motion of water in the presence of a symmetric cylinder floating freely (see
Fig.~1) be described by $\varphi$ even in $x$. Then equation \eqref{eq:5'} implies
that $z_1 = 0$, which means that such a potential does not comport with the free
sway motion of a symmetric cylinder. However, there are examples of motionless
symmetric cylinders with two immersed parts that trap modes at particular
frequencies; see \cite{4}, \S\S~5 and 6. 

Now, equality \eqref{eq:eq+} turns into
\begin{equation}
\nu \int_F |\varphi|^2 \, \D{}x = \int_W \big( |\nabla \varphi|^2 + k^2 |\varphi|^2
\big) \,\D{}x \D{}y \, , \label{eq:eq++=}
\end{equation}
and so the uniqueness theorem proved in \S~3 of \cite{4}, where the case of fixed
$B_-$ and $B_+$ was discussed, is applicable to the problem under consideration
here. The corresponding assertion formulated below is analogous to Proposition~3,
but involves property $\Omega_+$ instead of $\Omega_-$; namely, property $\Omega_+$
holds for $\omega$ if the inequalities
\[ \pi (2 m + 1) / 2 \leq \ell \, b \leq \pi (m + 1) , \ \mbox{where} \ \ell = (\nu^2 - 
k^2)^{1/2} \ \mbox{and} \ \nu = \omega^2 / g ,
\]
take place for some $m = 0,1,\dots$. Thus, these properties are of the same kind;
either of them is valid when $\ell \, b$ belongs to one interval in a sequence;
these sequences for $\Omega_+$ and $\Omega_-$ are complementary and their intervals
have common endpoints.

\begin{proposition}\label{propos.4}
Let $W$ have infinite depth and let the domain $\widehat{B}$ satisfy the assumptions
of Proposition 3. If property $\Omega_+$ holds for $\omega$, then problem
\eqref{eq:1}--\eqref{eq:5} and\/ \eqref{eq:6} has only a trivial solution $(\varphi,
\mathbf{z})$ for this $\omega$ provided $\varphi$ is even in $x$ and $\mathbf{z} =
(z_1, 0, 0)^\transp$.
\end{proposition}

Comparing Propositions 3 and 4 with the uniqueness theorem proved in \cite{4}, we
see that for a cylinder with two immersed parts the frequency intervals, where the
potential (either odd or even in $x$) is trivial, are the same irrespective whether
the cylinder floats freely or is fixed and depend only on the potential's parity.

\subsubsection{Heave motion}

In this case $\mathbf{z} = (0, z_2, 0)^\transp$, $z_2 \in \RR$, and so the boundary
condition \eqref{eq:4} takes the form
\begin{equation}
\pd{\mathbf{n}} \varphi = \omega N_2 z_2 \quad \mbox{on} \ S \, , \label{eq:4''}
\end{equation}
and system \eqref{eq:5} splits. Its second equation is as follows
\begin{equation}
\big( \omega^2 I^M - g I^D \big) z_2 = - \omega \int_{S} \varphi \, N_2 \, \D{}s \, ,
\label{eq:5''+}
\end{equation}
whereas the second and third equations turn into the orthogonality condi\-tions
\begin{equation}
\int_{S} \varphi \, N_1 \, \D{}s = 0 \, , \quad \int_{S} \varphi \, N_3 \, \D{}s = 0
\, . \label{eq:5''++}
\end{equation}
The last one is a consequence of the equality $I^D_x = 0$, which follows from the
symmetry about the $y$-axis; indeed, $x^{(0)} = 0$ in this case. Furthermore,
relation \eqref{eq:eq} takes the form:
\begin{equation}
\big( \omega^2 I^M - g I^D \big) z_2^2 = \nu \int_F |\varphi|^2 \, \D{}x - \int_W
\big( |\nabla \varphi|^2 + k^2 |\varphi|^2 \big) \,\D{}x \D{}y \, .
\label{eq:eq++}
\end{equation}

Assuming that $\varphi (x,y)$ is odd in $x$, we see that the right-hand side of
\eqref{eq:5''+} vanishes because the integrand has the same parity as $\varphi$.
Then the last equality turns into \eqref{eq:eq++=}, and so we are in a position to
refer to the considerations used in \cite{4}, \S~3. In this way, we arrive at the following assertion (it is
similar to Propositions~3 and 4) for the present problem.

\begin{proposition}\label{propos.5}
Let $W$ have infinite depth and let the domain $\widehat{B}$ satisfy the assumptions
of Proposition 3. If property $\Omega_-$ holds for $\omega$, then problem
\eqref{eq:1}--\eqref{eq:5} and\/ \eqref{eq:6} has only a trivial solution $(\varphi,
\mathbf{z})$ for this $\omega$ provided $\varphi$ is odd in $x$ and $\mathbf{z} =
(0, z_2, 0)^\transp$.
\end{proposition}

\noindent {\it Proof.} According to considerations in \cite{4}, \S~3, property
$\Omega_-$ guarantees that $\varphi$ (satisfying \eqref{eq:eq++=} and odd in $x$) is
trivial. Then the boundary condition \eqref{eq:4''} yields that $\mathbf{z}$ is also
trivial. $\square$

\vspace{1mm}

Let us turn to the case when $\varphi (x,y)$ is even in $x$ and apply the approach
used in the proof of Proposition~3 for determining the frequency intervals when
equality \eqref{eq:eq++} cannot hold if $\varphi$ is non-trivial. Considerations
based on property $\Omega_+$ show that the right-hand side of \eqref{eq:eq++} is
strictly negative when this property is valid. On the other hand, the left-hand side
of this equality is nonnegative if $\omega^2 \geq g I^D / I^M$; this yields the
following.

\begin{proposition}\label{propos.6}
Let $W$ have infinite depth and let the domain $\widehat{B}$ satisfy the assumptions
of Proposition 3. If property $\Omega_+$ holds for $\omega \geq \sqrt{g I^D / I^M}$,
then problem \eqref{eq:1}--\eqref{eq:5} and\/ \eqref{eq:6} has only a trivial
solution $(\varphi, \mathbf{z})$ for this $\omega$ provided $\varphi$ is even in $x$
and $\mathbf{z} = (0, z_2, 0)^\transp$.
\end{proposition}

The essential distinction between this assertion and Proposition 5 is as follows.
Along with property $\Omega_+$, the assumptions of Proposition 6 include the
inequality $\omega^2 \geq g I^D / I^M$, whereas no condition other than property
$\Omega_-$ is imposed in Proposition 5.

\subsubsection{Roll motion}

In this case $\mathbf{z} = (0, 0, z_3)^\transp$, $z_3 \in \RR$, and so the boundary
condition \eqref{eq:4} takes the form
\begin{equation}
\pd{\mathbf{n}} \varphi = \omega N_3 z_3 \quad \mbox{on} \ S \, , \label{eq:4''r}
\end{equation}
and system \eqref{eq:5} splits. The first and second equations turn into the
orthogonality conditions
\begin{equation}
\int_{S} \varphi \, N_1 \, \D{}s = 0 \, , \quad \int_{S} \varphi \, N_2 \, \D{}s = 0
\, . \label{eq:5''++r}
\end{equation}
Again, the last one is a consequence of the equality $I^D_x = 0$, which follows from
symmetry of $\widehat{B}$ about the $y$-axis. The third equation is as follows
\begin{equation}
\left[ \omega^2 I^M_2 - g \big( I^D_{xx} + I^B_y \big) \right] z_3 = - \omega
\int_{S} \varphi \, N_3 \, \D{}s \, , \label{eq:5''+r}
\end{equation}
where $N_3 (x, y) = x N_2 - (y - y^{(0)}) N_1$ and $N_3 (x, y) = - N_3 (-x, y)$ for
$(\pm x, y) \in S_\pm$ in view of symmetry. Now, relation \eqref{eq:eq} takes the
form:
\begin{equation}
\left[ \omega^2 I^M_2 - g \big( I^D_{xx} + I^B_y \big) \right] z_3^2 = \nu \int_F
|\varphi|^2 \, \D{}x - \int_W \big( |\nabla \varphi|^2 + k^2 |\varphi|^2 \big)
\,\D{}x \D{}y \, .
\label{eq:eq++r}
\end{equation}

Assuming that $\varphi (x,y)$ is even in $x$, we see that the right-hand side of
\eqref{eq:5''+r} vanishes (cf. \S~4.2.2), and the last equality turns into
\eqref{eq:eq++=}. Again we are in a position to refer to the considerations used in
\cite{4}, \S~3, which lead to.

\begin{proposition}\label{propos.7}
Let $W$ have infinite depth and let the domain $\widehat{B}$ satisfy the assumptions
of Proposition 3. If property $\Omega_+$ holds for $\omega$, then problem
\eqref{eq:1}--\eqref{eq:5} and\/ \eqref{eq:6} has only a trivial solution $(\varphi,
\mathbf{z})$ for this $\omega$ provided $\varphi$ is even in $x$ and $\mathbf{z} =
(0, 0, z_3)^\transp$.
\end{proposition}

Comparing this assertion and Proposition 4, we see that for a freely floating
symmetric cylinder with two immersed parts property $\Omega_+$ guarantees uniqueness
at all frequencies satisfying it provided $\varphi$ is even in $x$ and the cylinder
executes either sway or roll motion.

In the case when $\varphi (x,y)$ is odd in $x$, we again apply the approach used for
the proof of Proposition~6. Since  $I^D_x = 0$ in view of symmetry of $\widehat{B}$,
the assumption that $\mathbf{K}'$ is a positive definite matrix (see \S~2) yields
the inequality $I^D_{xx} + I^B_y > 0$. Thus we obtain.

\begin{proposition}\label{propos.8}
Let $W$ have infinite depth and let the domain $\widehat{B}$ satisfy the assumptions
of Proposition 3. If property $\Omega_-$ holds for
\begin{equation}
\omega \geq \sqrt{g \big( I^D_{xx} + I^B_y \big) / I^M_2} \, ,
\label{eq:last}
\end{equation}
then problem \eqref{eq:1}--\eqref{eq:5} and\/ \eqref{eq:6} has only a trivial
solution $(\varphi, \mathbf{z})$ for this $\omega$ provided $\varphi$ is odd in $x$
and $\mathbf{z} = (0, 0, z_3)^\transp$.
\end{proposition}

As in the case of Propositions 5 and 6, the distinction between Propositions 7 and 8
is that the latter one includes the extra inequality \eqref{eq:last}.

In conclusion of this section, it should be mentioned that there are analogues of
Propositions 3--8 for water of finite depth; their formulations are similar to those
above.

\vspace{-2mm}

\section{Discussion}

First, the obtained results guarantee uniqueness of a solution in proper classes of
functions for the problem of scattering of obliquely incoming plane waves by an
infinitely long cylinder floating freely. This problem is left for future research.

Second, equations \eqref{eq:5''+} and \eqref{eq:5''+r} demonstrate that one may
expect the existence of trapped modes for some particular bodies in heave and roll
motion. To find these bodies and the corresponding trapped modes is another
interesting problem for future research. A hint to its solution (at least for $k=0$)
can be found in the article \cite{KM2}.


\vspace{-1mm}

\begin{thebibliography}{99}

\bibitem{NGK10} N. Kuznetsov, \textit{On the problem of time-harmonic water waves
in the presence of a freely-floating structure.} --- St Petersburg Math. J. {\bf 22}
(2011), 985--995.

\bibitem{John1} F. John, \textit{On the motion of floating bodies, I.} --- Comm.
Pure Appl. Math. {\bf 2} (1949), 13--57.

\bibitem{U} F. Ursell, \textit{Some unsolved and unfinished problems in the theory
of waves.} --- Wave Asymptotics, eds. P.A. Martin, G.R. Wickham, Cambridge,
Cambridge University Press (1992), 220--244.

\bibitem{CK} D. Colton, R. Kress, \textit{Integral Equation Methods in Scattering
Theory.} --- Philadelphia, SIAM (2013).

\bibitem{KM} N. Kuznetsov, O. Motygin, \textit{On the coupled time-harmonic motion
of water and a body freely floating in it.} --- J. Fluid Mech. {\bf 679} (2011),
616--627.

\bibitem{KM1} N. Kuznetsov, O. Motygin, \textit{On the coupled time-harmonic motion
of deep water and a freely floating body: trapped modes and uniqueness theorems.}
--- J. Fluid Mech. {\bf 703} (2012), 142--162.

\bibitem{KM2} N. Kuznetsov, O. Motygin, \textit{Freely floating structures trapping
time-harmonic water waves.} --- Quart. J. Mech. Appl. Math. {\bf 68} (2015),
173--193.

\bibitem{John} F. John, \textit{On the motion of floating bodies, II.} --- Comm.
Pure Appl. Math. {\bf 3} (1950), 45--101.

\bibitem{KMV} N. Kuznetsov, V. Maz'ya, B. Vainberg, \textit{Linear Water Waves: A
Mathematical Approach.} --- Cambridge, Cambridge University Press (2002).

\bibitem{MMI} M. McIver, \textit{An example of non-uniqueness in the
two-dimensional linear water-wave problem.} --- J. Fluid Mech. {\bf 315} (1996),
257--266.

\bibitem{4} N. Kuznetsov, R. Porter, D.\,V. Evans, M.\,J. Simon, \textit{Uniqueness
and trapped modes for surface-piercing cylinders in oblique waves.} --- J. Fluid
Mech. {\bf 365} (1998), 351--368.

\bibitem{K} N. Kuznetsov, \textit{Uniqueness in the problem of an obstacle in
oblique waves.} --- C. R. Mecanique. {\bf 331} (2003), 183--188.

\bibitem{KM3} N. Kuznetsov, O. Motygin, \textit{On the coupled time-harmonic motion
of water and a body freely floating in it.} --- J. Fluid Mech. {\bf 679} (2011),
616--627.

\end{thebibliography}
\end{document}